\documentclass[12pt]{article}

\makeatletter
\def\@cite#1#2{\nolinebreak$^{[\scriptstyle #1\if@tempswa , #2\fi]}$}
\def\@citex[#1]#2{\if@filesw\immediate\write\@auxout{\string\citation{#2}}\fi
  \def\@citea{}\@cite{\@for\@citeb:=#2\do
    {\@citea\def\@citea{,\penalty\@m}\@ifundefined
       {b@\@citeb}{{\bf ?}\@warning
       {Citation `\@citeb' on page \thepage \space undefined}}%
{\csname b@\@citeb\endcsname}}}{#1}}

\@definecounter{footnote}

\begin{document}

\newcommand{\be}{\begin{equation}}
\newcommand{\ee}{\end{equation}}

\parskip 15pt
\baselineskip 18pt

\centerline{\bf\large Lanchester combat models}

\vskip 0.2in \centerline{\large Niall MacKay}

\centerline{Department of Mathematics, University of York}

\centerline{May 2005}

\vskip 0.2in

\noindent{\bf 1. The Lanchester aimed-fire model}

\noindent As a simple, deterministic model\cite{Lanch} of a
battle, suppose that $R(t)$ red and $G(t)$ green units\footnote{We
shall worry later on about what exactly constitutes a `unit'.}
begin fighting at $t=0$, and that each unit destroys $r$ or $g$
(the {\bf fighting effectiveness}) enemy units in one unit of
time, so that \be\label{LDFM} {dR\over dt} = -gG\;,\hspace{0.6in}
{dG\over dt} = -rR\;.\ee Rather than solve directly, we eliminate
the explicit $t$-dependence by dividing the second equation by the
first and then separating variables: then \be \int rR\,dR = \int
gG\,dG\;,\ee and we see that \be rR^2-gG^2 = {\rm
constant.}\label{DFS}\ee This has remarkable implications. Since
this quantity never changes sign, only one of $R$ and $G$ can ever
be zero. If the initial value of (\ref{DFS}) is positive, for
example, only $G$ can equal zero, and so only red can
win\footnote{To show that red actually does win, we would need to
solve the equations explicitly for $R(t)$ and $G(t)$, which is
most easily done by differentiating one of (\ref{LDFM}),
substituting into the other, and solving the resulting second
order ODE.}. So this quantity is key in determining who will win
the battle, and a side's fighting strength -- we shall call $rR^2$
and $gG^2$ the red and green forces' respective {\bf fighting
strengths} -- varies as the units' fighting effectiveness times
the {\em square} of their numbers.

Consider an example. Suppose red begins with twice as many units
as green, $R_0=2G_0$, but the green units are three times as
effective, $g=3r$. Then \be rR^2-gG^2 = r(2G_0)^2-3rG_0^2= rG_0^2
>0\;,\ee and (perhaps rather counter-intuitively) the reds win.
The model assumes that the battle is a homogeneous mixing of
forces, or, in other words, a constant bloody {\em m\^el\'ee}, and
in such engagements the tactical conclusion is simple: if your
strength is in numbers, then you need to fight in this way,
bringing all your units to engage with the enemy's as rapidly as
possible. Conversely, if your units are fewer in number but more
effective, then you need the tactics which will prevent a {\em
m\^el\'ee}, allowing you to pick-off opponents, and preventing the
enemy from bringing all his units to bear ( -- like Hannibal at
Cannae).

This last point becomes clearer if we consider what would have
happened had green been able to divide the reds into two equal
forces and engage them sequentially. At the end of the first
engagement, between $G_0$ greens and $R_0=G_0$ reds, $G_1$ greens
remain, where \be rG_0^2 - 3r G_0^2= -3r G_1^2
\hspace{0.5in}\Rightarrow \hspace{0.5in} G_1 = \sqrt{2\over 3}
G_0\;,\ee and then in the second engagement \be r G_0^2 - 3r
{2\over 3}G_0^2 = -rG_0^2<0\;. \ee Green has now won, with
$3rG_2^2=rG_0^2$ and thus $\sqrt{1\over 3}$ or nearly 60\% of its
original forces left when all the reds have been destroyed -- an
amazing turnaround. This becomes even more striking with an
$N$-fold (rather than two-fold) division of red forces: after the
simple resulting iteration, green now wins with a final number
$G_F$ of units remaining, where \be gG_F^2 = gG_0^2 - {1\over
N}rR_0^2\;,\ee so that the $N$-fold division has reduced red's
fighting strength $N$-fold. Again, this is a classic military
maxim: you should (almost) never divide your forces\footnote{It is
perhaps worth noting that if red and green both divided their
forces, and fought two separate, simultaneous, identical battles,
the outcome would necessarily be the same as for our single,
original battle. The argument is rather like Galileo's for the
falling cannonballs: making an imaginary or a real division of the
battle makes no difference to its course or outcome.}.

 \noindent{\bf 2. Background}

\noindent The above is known as the Lanchester aimed- (or
directed-) fire combat model, and it seems to me to exemplify all
that's best in a simple mathematical model. Most modern warfare
simulations are, of course, stochastic, heterogeneous and complex;
and will, {\em if} their myriad assumptions are correct, give much
better predictions. But the Lanchester model, in
contrast\cite{Lucas3}, has the virtues of simplicity: it makes
strong simplifying assumptions, which nevertheless are (at least
sometimes) close to being realizable, and the model brings out,
through a subtle process, some stark conclusions.

So how close a fit is the model to past battles? Lanchester
originally applied it to Nelson's tactics at Trafalgar, and
intended that it should describe aerial combat, but for detailed
fits of data one has to look to battles such as Iwo
Jima\cite{Engel}, the Ardennes\cite{Fricker} and
Kursk\cite{Lucas1}. (Indeed my values and notation above are
intended approximately to match Kursk, the great Russian-German
tank battle of 1943.) In fact the fit seems not to be so
good\cite{Lucas1}, but my impression after exploring the
literature is that this is to be expected: the conditions implied
by the model, of constant tactically-blind slaughter, are hardly
those which allow the collection of accurate time-series data; and
conversely where good data exist they are more likely to describe
a series of different small engagements. Those who teach military
tactics, however, still value Lanchester's model and its
generalizations\cite{Lucas3}, because, above all, they provoke
careful thought about the consequences of the conditions of
engagement.

It seems to me that the model offers an excellent pedagogical tool
for curriculum enrichment at least down to A-level (where simple
separable first-order ODEs are in the core curriculum). Moreover,
it is easily extended to suit students with no calculus (see
section 3.1 below). Perhaps its main virtue is that so many of the
natural questions one can ask about it have
mathematically-tractable answers. We investigate some of these in
the next section.

\noindent{\bf 3. Further explorations}

\noindent{\em 3.1. A calculus-free version}

\noindent A model in the form of discrete recursion relations not
only avoids calculus but is also essential in modelling battles
between small numbers of units, or battles effectively fought as a
sequence of discrete engagements (like salvos in a naval
battle\cite{Hughes}). The difference equations \be
R_{n+1}=R_n-gG_n\;,\hspace{0.6in} G_{n+1}=G_n-rR_n \;,\ee whose
continuous-time limit is (\ref{LDFM}), lead to \be
rR_{n+1}^2-gG_{n+1}^2 = (1-rg)\left(rR_n^2-gG_n^2\right)\,.\ee But
$rg\ll 1$, and indeed $(1-rg)^T\simeq 1$ for the duration $T$ of
the battle (which is $\sim 1/r$ or $1/g$). Thus (\ref{DFS})
remains approximately conserved, and the same conclusions apply.

\noindent{\em 3.2. The Lanchester unaimed-fire model}

\noindent There are various scenarios in which a side's fighting
strength is proportional to its numbers rather than their square.
Indeed any model of the form \be {dR\over dt} =
-F_G\;,\hspace{0.6in} {dG\over dt} = -F_R\;,\label{LUFM}\ee where
$F_G/F_R$ is independent of $R$ and $G$, gives this result.
Lanchester considered that $F=$constant modeled (`ancient')
hand-to-hand combat without firearms, while a modern model with
this property would be that of artillery fire or bombardment, in
which, if $R$ red guns fire with effectiveness $r$ at random into
an area $A_G$ in which there are known to be $G$ green units
(which therefore have density $G/A_G$ per unit area), we have
$F_R=-rRG/A_G$. Either way, it is now $ \rho R-\gamma G $ (for
some constants $\rho$ and $\gamma$) which is constant, and the
fighting strength of unaimed-fire units is proportional to their
numbers and effectiveness but inversely proportional to the
enemy's area of dispersal. Once again, there are some simple
implications. For example, in a tank battle you may wish to engage
the enemy with aimed-fire from well-dispersed and -disguised
positions, forcing the enemy to use unaimed fire. Thus you might
emplace your guns or tanks, perhaps within the edge of a wood.
However, if the enemy later becomes able to aim his fire, he may
be able to divide your forces; you will need to be able to regroup
quickly.

\noindent{\em 3.3. Mixed forces}

\noindent Of course fighting forces have always been composed of
units of different effectiveness in different numbers. Suppose the
green forces are composed of two types, so $G=G_1+G_2$, with
effectiveness $g_1$ and $g_2$. Then \be {dR\over dt} =
-g_1G_1-g_2G_2\;,\hspace{0.4in} {dG_1\over dt} = -rR{G_1\over
G}\,, \hspace{0.3in}{dG_2\over dt} = -rR{G_2\over G}
\;.\label{mixed}\ee Although this looks more complicated, in fact
the outcome is surprisingly simple: one can check straighforwardly
that \be {d\over dt}\left(
rR^2-(g_1G_1+g_2G_2)(G_1+G_2)\right)=0\,,\ee so that the conserved
quantity analogous to (\ref{DFS}) is now \be rR^2- g_{\rm ave}
G^2\;,\ee where \be g_{\rm ave}= {g_1G_1+g_2G_2\over G}\ee is the
average effectiveness of the units.
 So to calculate the fighting strength of mixed units
one should simply use the number of autonomous units with their
average effectiveness (in contradiction with Lanchester's original
claims, as has been pointed out by Lepingwell\cite{Lepingwell},
whose article is an excellent introduction to Lanchester models).

\noindent{\em 3.4. The meaning of a `unit'}

\noindent This leads us to an interesting controversy in the
military use of the model\cite{Homer}. Sometimes the use of
Lanchester-type models is taught in which units' {\em numbers} are
first weighted with a `fighting power' -- for instance, 1000
troops together with 20 tanks, the latter multiplied by a
`fighting power' of, say, 30, to give a total of 1600 units. As
others have pointed out, and as we can see from the above, this
cannot be correct -- the Lanchester model's main point is to
distinguish the importance of numbers from that of fighting
effectiveness. But it does lead us to a natural and fundamental
question, which should by now be worrying the reader: what exactly
do we mean by a `unit'?

For example, consider an infantry engagement. Is the basic unit
the individual soldier, or perhaps the section (about 10 men),
platoon (3-4 sections), or company (3-4 platoons)? Does this
choice change if (for example) each section fights with its own
armoured personnel carrier (APC)? Above all, does this choice
affect the outcome of the battle in our model? In fact it does
not. For suppose we take green's basic unit to be $N$ troops. Then
its numbers are scaled by $1/N^2$ and its effectiveness by $N$.
But red's effectiveness against green units is now scaled by
$1/N$, and the whole of the difference in fighting strengths
(\ref{DFS}) scales by $1/N$. Its sign, and the course and outcome
of the battle, do not change.

But $N$ troops may gain some advantage by fighting as a group. It
may be offensive, in which case their effectiveness scales better
than $N$-fold, or defensive, in which case their opponents'
effectiveness will scale worse than $1/N$ (for example, in
hand-to-hand fighting, a small group fighting in a cluster and
thus defending each other's backs). But of course we now rapidly
run into the limitations of the model, for battles with different
types of units are rarely homogeneous -- indeed it is a tactical
imperative to use units so as to maximise their fighting
advantages. A section fighting in an APC may be better protected
against other infantry, but it is also now participating in a
different, armoured battle.

\noindent{\em 3.5. Support troops}

\noindent During the battle for France in 1940, Churchill visited
the French HQ and, while being shown on a map the location of the
German breakthrough, asked {\em `Ou est la masse de manoeuvre?'}
(`Where is the strategic reserve?'). That there was none -- {\em
`aucune'} -- confirmed to him that the battle was lost. Tactical
considerations, beyond the scope of Lanchester models, suggest
that one should always maintain a reserve. But support and/or
reserve troops may also increase the efficiency of fighting
troops, and we must ask which is correct: to aggregate support
troops within the overall numbers, or to deal with them
separately. A thorough analysis clearly requires the latter.

Suppose that we split our $N$ units into $P$ fighting and $N-P$
support units. We do so because support units increase the
effectiveness of fighting units, so let us assume that this is
proportional to the support ratio ${N-P\over P}$, so that the $P$
fighting units' effectiveness $f =f_0 {N-P\over P}$. Then the
fighting strength $fP^2$ is maximized by $P=N/2$ -- that is, half
of our units should be in support roles. More generally, if we
take $f = f_0 \left({N-P\over P}\right)^\kappa$, then the fighting
strength is maximized by $P=(1-\kappa/2)N$. For example, a ratio
of one headquarters to three fighting companies, or $P/N=3/4$ and
$\kappa=1/2$, might correspond to an implicit belief that a 10\%
increase in such support per fighting unit leads to an approximate
5\% increase in its effectiveness.

In the Ardennes campaign\cite{Bracken}, the Allies used a much
greater support ratio than the Germans (about 0.8 as opposed to
0.5), but the fighting effectiveness of their total numbers -- in
our model above, $f_{\rm
tot}=f_0(\kappa/2)^\kappa(1-\kappa/2)^{2-\kappa}$ -- turns out to
have been about the same as the Germans'. Of course, this does not
mean that the Allies should have reduced their support ratio. If
we assume that this was optimal (perhaps because of more plentiful
{\em mat\'eriel}), so that $\kappa\simeq 0.89$, then to have
reduced ${N-P\over P}$ to that of the Germans would have reduced
their $f_{\rm tot}$ by about 2\%. It is not clear whether, in the
model, this would have been within the margin of victory!

 \noindent{\em 3.6.
Bracken's generalized model}

\noindent As we discussed in section two, attempts in the
literature to fit either the aimed- or the unaimed-fire model to
data have only been partly successful. If we are willing to move
away from basing our model on clear assumptions (always a
dangerous path), we can instead use a generalized model whose
parameters we then fit to the data. In this spirit Bracken
proposed\cite{Bracken} \be {dR\over dt} = -g
R^qG^p\;,\hspace{0.6in} {dG\over dt} =
-rR^pG^q\;,\label{Bracken}\ee for $p$ and $q$ to be empirically
determined.  The conserved quantity (by eliminating $t$,
separating and integrating) is \be gG^\alpha - rR^\alpha\,,\ee
where $\alpha=1+p-q$. (The Lanchester aimed-fire model
(\ref{LDFM},\ref{DFS}) corresponds to $p=1,q=0$ and thus to
$\alpha=2$, the unaimed-fire model to $p=q=1$ and thus
$\alpha=0$.)

Clearly if $\alpha>1$ then numbers matter more than effectiveness,
and conversely. It might be natural to assume that $p>q$ (and thus
$\alpha>1$), so that one's rate of loss scales more quickly with
the enemy's numbers than with one's own, but in fact an enormous
variety of best-fitting $p,q$ have been found\cite{Fricker,Lucas1}
by looking at various battles and with differing measures of
numbers and effectiveness -- which certainly indicates that this
type of model is of little use in advance prediction of a
particular battle's outcome! Hartley's analysis\cite{Hartley} (as
reported by Lucas and Dinges\cite{Lucas2}) of a range of battles
suggests $p=0.45$, $q=0.75$.

\noindent{\em 3.7. Asymmetric warfare}

\noindent This result, that one's own numbers can have at least as
strong an effect on one's casualties as do enemy numbers, might
lead us to add exponential decay terms in both equations of the
original model. Instead, let us  revisit a fundamental assumption
in all of the foregoing: that the correct model is symmetric under
$R\leftrightarrow G$ and their associated parameters.

Suppose red has huge numerical superiority and
 green has vastly greater effectiveness -- that is, $g\gg r$ but
$G_0\ll R_0$. If, further, $rR_0^2 \gg gG_0^2$ then the Lanchester
aimed-fire model would give a clear win for red. However, in the
early stages of the battle there will be many reds for each green,
which will find it easy to acquire targets, and thus have a kill
rate proportional to both $R$ and $G$. The reds, meanwhile, will
kill in proportion to $G$ but will not be able to make their
numbers tell. So let us take \be\label{asymmetric} {dR\over dt} =
-gGR/R_0\;,\hspace{0.6in} {dG\over dt} = -rG\;.\ee Then $G$
declines slowly and exponentially, $G=G_0e^{-rt}$, and we can now
solve for $R$ to give \be R=R_0 \exp\left\{ {gG_0\over
rR_0}(e^{-rt}-1)\right\}\,.\ee Now consider the situation when
$e^{-rt}={1\over 2}$. The greens have been reduced to half their
original numbers, $G={1\over 2} G_0$, while
$R=R_0\exp\left(-{1\over 2} {gG_0\over rR_0}\right)$. So, if
$gG_0\gg rR_0$, -- in other words, if green would clearly win in
the unaimed-fire, linear-law model, -- then it will win at least
the early stages of this battle\footnote{Of course one then has to
switch to a standard aimed-fire model for the later stages. Green
may still lose the battle in the end!}, even though it would lose
under a conventional aimed-fire, square-law model. Thus, with
technologically-superior forces in a `target-rich environment',
even under homogeneous, aimed-fire conditions, it may be a linear-
rather than a square-law battle which is being fought.

\noindent{\bf 4. Concluding remarks}

\noindent Of course there have been many generalizations of the
Lanchester models, but most of these lack the simplicity and force
of the original. One can include exponential or linear growth or
attrition\cite{Taylor} in the original equations, corresponding to
defence\cite{Hughes}, friendly fire or reinforcement, and all of
these are easily soluble with A-level or first-year undergraduate
techniques. A model which mixes aimed with unaimed fire would,
however, be nonlinear (a simple Lotka-Volterra model), and perhaps
not for pre-university students. Further, any of these properties
can appear asymmetrically in the coupled equations. There are many
websites dealing with various of the possibilities.

Lanchester models provide an excellent example of the strengths
(and weaknesses) of simple mathematical modelling. Further, as we
have seen, the basic model leads to many other `What if...?'
questions which can be easily investigated. Many more such
questions can be asked, and of course once one begins numerical
simulations the possibilities are endless. To model warfare can
seem more politically challenging (not to say incorrect) than to
use ecology or epidemiology, but a little understanding of how
military planners arrive at their tactical conclusions can also
strip away mystique and (for this author, at least) expose some of
the subject's limitations!

{\bf Acknowledgments}. My original exposure to Lanchester models
was through a talk by Sir John Kingman, and I should like to thank
him for stimulating my interest. Thanks to Simon Eveson and
Charles Young for comments on the manuscript.

\end{document}